\documentclass[12pt]{article} \usepackage{amssymb,epsfig,psfrag,amsmath,art,cite} \usepackage{multicol}
\usepackage{subfigure}[1995/03/06] \newfont{\smalll}{cmr8}

\def\IR{\mathbb{R}}
\def\IS{\mathbb{S}}

\def\span{\mathrm{span}}

\def\IS{\hbox{I\hskip-.1em S}}
\def\IC{\hbox{C\hskip-
.5em\raise.5ex\hbox{$\scriptscriptstyle\mid$}}\ }
\def\Ic{\hbox{\smalll C\hskip-
.5em\raise.3ex\hbox{$\scriptscriptstyle\mid$}}\ }

\def\T={\buildrel {\scriptscriptstyle\triangle} \over =}

\def\sqr#1#2{{\vcenter{\vbox{\hrule height.#2pt\hbox{\vrule
width.#2pt height#1pt \kern#1pt\vrule width.#2pt}\hrule
height.#2pt}}}}
\def\square{\mathchoice\sqr64\sqr64\sqr33\sqr33}

\def\block-diag{\mathop{\rm block{\scriptstyle -}diag}}

\def\pmbb#1{\setbox0=\hbox{#1}\raise 0.5ex\box0}

\def\norm#1{\|#1\|}

\newcommand{\bequ}{\begin{eqnarray}}
\newcommand{\eequ}{\end{eqnarray}}
\newcommand{\mT}{^\mathrm{T}}

\newcommand{\rom}{\mathrm}

\newcommand {\teq}      {\triangleq}

\newcommand {\beq}      {\begin{equation}}
\newcommand {\eeq}      {\end{equation}}

\newcommand{\Proj}{{\mathrm{Proj}}}

\def\IR{{\mathbb R}}
\def\IC{{\mathbb C}}

\def\IS{{\mathbb S}}

\newcommand{\union}{\mathop{\cup}}

\renewcommand\theenumi{\roman{enumi}}
\renewcommand\theenumii{\alph{enumii}}
\renewcommand\theenumiii{\arabic{enumiii}}
\renewcommand\theenumiv{\Alph{enumiv}}

\def\sc{s_{\rm c}}

\def\sc{s_{\rm c}}

\begin{document}

\title{\Large{Improving Transient Performance of Adaptive Control\\Architectures using Frequency-Limited\\System Error Dynamics}\\}
\def\YucelenJohnson
{\begin{tabular}{c}
   Tansel Yucelen \\
   Department of Mechanical and Aerospace Engineering \\
   Missouri University of Science and Technology \\
   400 W. 13th St., Rolla, MO 65401 \\
   PHONE: (573) 341-7702 \\
   FAX: (573) 341-6899 \\
   yucelen@mst.edu \\
\\
\\
 \end{tabular}
\ \\[-1.5em]
\begin{tabular}{ccc}
Gerardo De La Torre & & Eric N. Johnson\\
School of Aerospace Engineering&  & School of Aerospace Engineering\\
Georgia Institute of Technology &  & Georgia Institute of Technology\\
270 Ferst Dr., Atlanta, GA 30332 & & 270 Ferst Dr., Atlanta, GA 30332\\
PHONE: (404) 385-4940 & & PHONE: (404) 385-2519 \\
FAX: (404) 894-2760 & & FAX: (404) 894-2760\\
glt3@gatech.edu& &eric.johnson@ae.gatech.edu
\end{tabular}
}
\author{\YucelenJohnson}\date{ October 2013 \vspace{0.75em} } \maketitle \baselineskip 12pt

\vspace{2.5em} \begin{center} {\bf Abstract} \end{center}

We develop an adaptive control architecture to achieve stabilization and command following 
of uncertain dynamical systems with improved transient performance. 
Our framework consists of a new reference system and an adaptive controller. 
The proposed reference system captures a desired closed-loop dynamical system behavior modified by a mismatch term 
representing the high-frequency content between the uncertain dynamical system and this reference system, 
i.e., the system error. 
In particular, this mismatch term allows to limit the frequency content of the system error dynamics, 
which is used to drive the adaptive controller. 
It is shown that this key feature of our framework yields fast adaptation without incurring 
high-frequency oscillations in the transient performance. 
We further show the effects of design parameters on the system performance, 
analyze closeness of the uncertain dynamical system to the unmodified (ideal) reference system, 
discuss robustness of the proposed approach with respect to time-varying uncertainties and disturbances, 
and make connections to gradient minimization and classical control theory. 


\vspace{1.5em}




\clearpage \baselineskip=24pt \setcounter{page}{1}


\section{Introduction}

\begin{figure}[b!] \center \epsfig{file=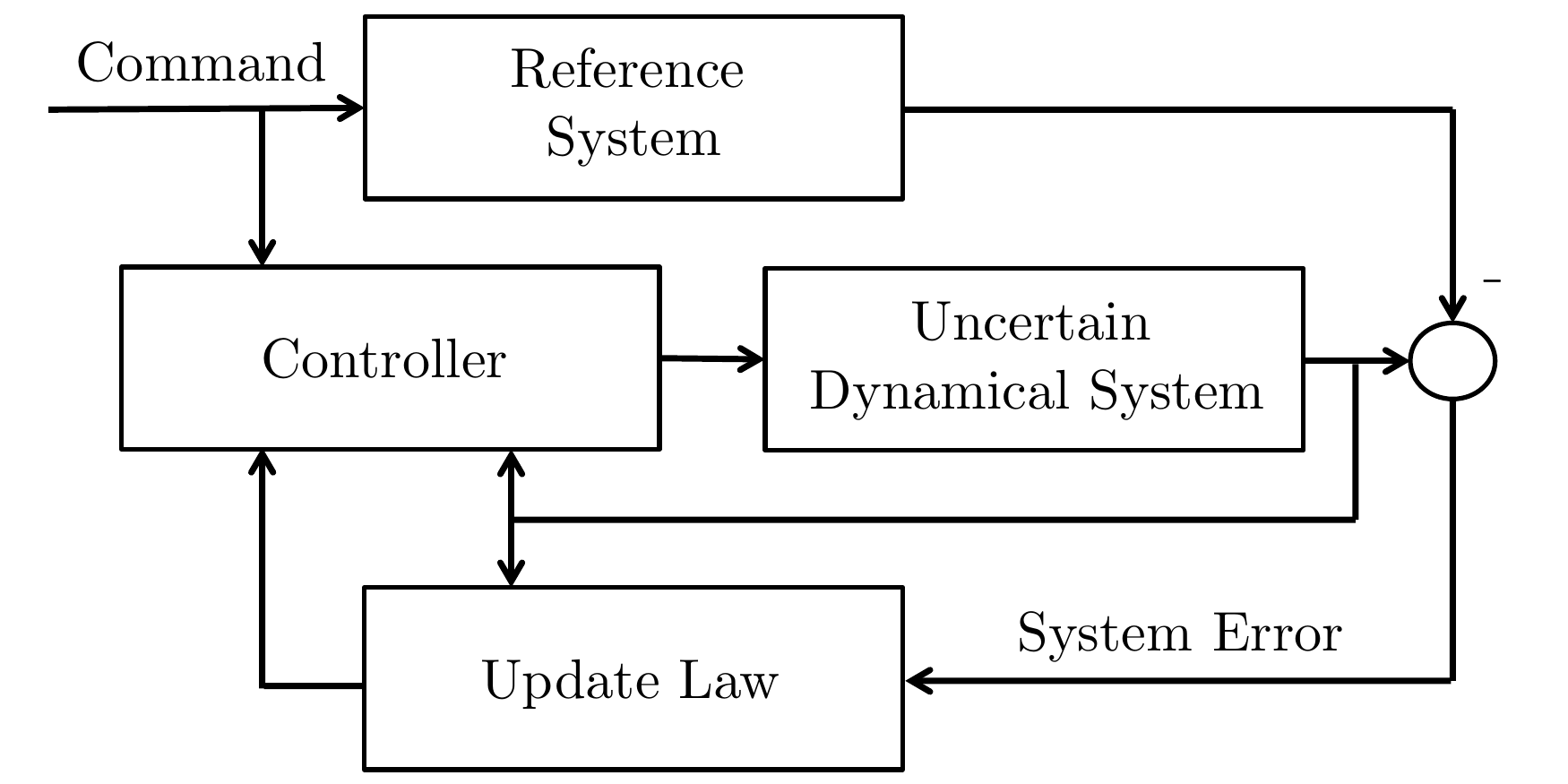, scale=0.65}
\caption{Block diagram of the model reference adaptive control scheme.}
\label{early:figure:1} \end{figure}

Although fixed-gain robust control design approaches can deal with uncertain dynamical systems, they require the knowledge of uncertainty bounds. 
Characterization of these bounds is not trivial due to practical constraints, because it requires extensive and costly verification and validation procedures. 
Furthermore, in the face of high uncertainty levels, system faults, or structural damage, these approaches may fail to satisfy a given system stabilization or command following requirement. 
On the other hand, adaptive control design approaches can effectively deal with these sources of uncertainties and require less modeling information than do fixed-gain robust control approaches. 
These facts make adaptive control theory a candidate for many applications.
The control framework of this work builds on a well-known and important class of adaptive controllers, specifically, model reference adaptive controllers. 
Whitaker \textit{et al.} [\citen{whitaker:1}, \citen{whitaker:2}] originally proposed the model reference adaptive control concept. 
In particular, model reference adaptive control schemes have three major components, namely, a reference system (model), an update law, and a controller (Figure \ref{early:figure:1}).
The reference system, in the classical sense, captures a desired closed-loop dynamical system behavior for which its output (resp., state) is compared with the output (resp., state) of the uncertain dynamical system. 
This comparison results in a system error signal used to drive the update law online. 
Then, the controller adapts feedback gains to minimize this error signal using the information received from the update law. 

From a practical standpoint, the output (resp., state) of the uncertain dynamical system can be far different from the output (resp., state) of the reference system during transient time (learning phase), although a model reference adaptive control scheme can guarantee that the distance between the uncertain dynamical system and the reference system vanishes asymptotically. 
This problem, so-called poor transient performance phenomenon, can be solved by increasing the learning rate of the update law, and hence, fast adaptation can be achieved in order to suppress uncertainties rapidly during transient time. 
However, update laws with high learning rates may yield to signals with high-frequency content, which can violate actuator rate saturation constraints and/or excite unmodeled system dynamics [\citen{oscil:2}, \citen{oscil:3}] resulting in system instability for practical applications. 
Hence, a critical trade-off between system stability and control adaptation rate exists in most adaptive control approaches, with some notable exceptions [\citen{notable:1}--\citen{notable:4}].
Specifically, the authors in \cite{notable:1} use a low-pass filter that subverts high-frequency oscillations attributable to fast adaptation, and their approach has guaranteed transient and steady-state performance. 
The authors in \cite{notable:2} present an adaptive control approach based on high learning rates that allows fast adaptation without hindering system robustness. 
Even though the methodologies documented in [\citen{notable:1}, \citen{notable:2}] are promising, they require the knowledge of a conservative upper bound on the unknown constant gain appearing in their uncertainty parameterization. 
While this conservative upper bound can be available for some applications, the actual upper bound may change and exceed its conservative estimate, for example, when an aircraft in a flight control application undergoes a sudden change in dynamics. 
In such circumstances, the performance of these adaptive controllers may be poor, because tuning them online with a new conservative upper bound is not possible.  
The author in \cite{notable:3} presents a modification to the reference system in order to solve the poor transient performance phenomenon, where a detailed analysis of this approach is given in \cite{notable:4}. 
Specifically, this modification is constructed by using a modification gain multiplied by the system error that is between the uncertain dynamical system and the modified reference system. 
In the limit as this modification gain goes to infinity, it is shown that the system error goes to zero in transient time. 
This approach can be used to effectively suppress uncertainties, however, for example, in the presence of exogenous low-frequency persistent disturbances, the transient performance of this approach may not be sufficient.  
Because, this disturbance may not be \textit{visible} to the update law, since the system error is (sufficiently) small due to a (sufficiently) large modification gain.

This work develops an adaptive control architecture to achieve stabilization and command following of uncertain dynamical systems with improved transient performance. 
The contribution of our framework comes from using a new reference system with an adaptive controller. 
The proposed reference system captures a desired closed-loop dynamical system behavior modified by a mismatch term
\begin{figure}[t!]   \center   \epsfig{file=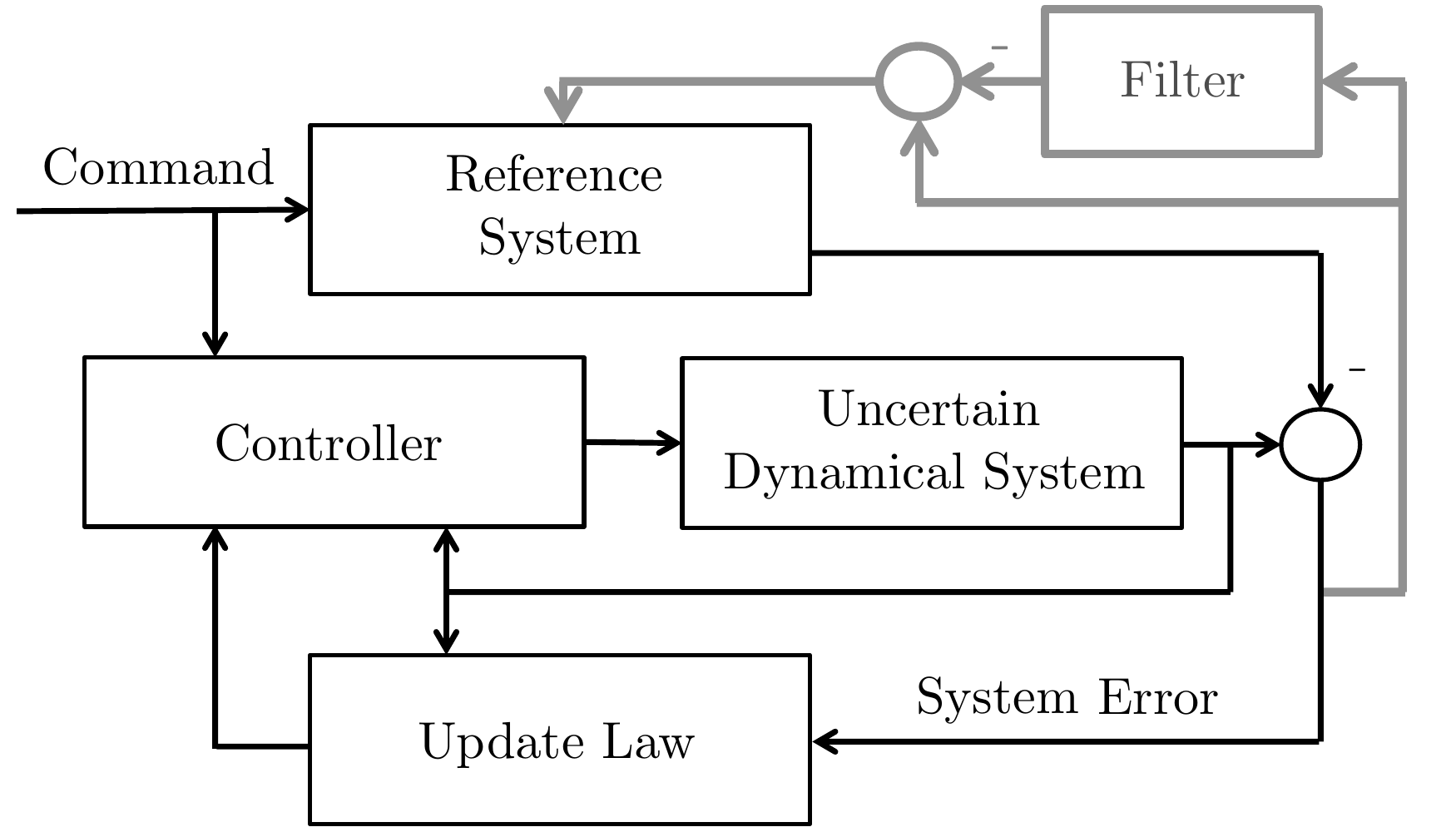, scale=0.65}
\caption{Block diagram of the proposed scheme. Note that the reference system is driven not only by the command but also by the difference between the system error and its (low-pass) filtered form representing the high-frequency content of the system error.}
\label{early:figure:2} \end{figure}
representing the high-frequency content between the uncertain dynamical system and this reference system, i.e., the system error (Figure \ref{early:figure:2}).
In particular, this mismatch term allows to limit the frequency content of the system error dynamics, which is used to drive the adaptive controller. 
That is, the update law learns through the low-frequency content of the system error, which constitutes a distinction over the approach in \cite{notable:3}. 
It is shown that this key feature of our framework yields fast adaptation without incurring 
high-frequency oscillations in the transient performance. 
We further show the effects of design parameters on the system performance, 
analyze closeness of the uncertain dynamical system to the unmodified (ideal) reference system, 
discuss robustness of the proposed approach with respect to time-varying uncertainties and disturbances, 
and make connections to gradient minimization and classical control theory. 


\ \section{Notation}

The notation used in this paper is fairly standard.
Specifically, $\IR$ denotes the set of real numbers,
$\IR^n$ denotes the set of $n \times 1$ real column vectors,
$\IR^{n \times m}$ denotes the set of $n \times m$ real matrices,
$\IR_+$ (resp., $\overline{\IR}_+$) denotes the set of positive (resp., nonnegative-definite) real numbers,
$\IR_+^{n \times n}$ (resp., $\overline{\IR}_+^{\hspace{0.1em} n \times n}$) denotes the set of $n \times n$ positive-definite (resp., nonnegative-definite) real matrices,
$\IS^{n \times n}$ denotes the set of $n \times n$ symmetric real matrices,
$\mathbb{D}^{n \times n}$ denotes the $n \times n$ real matrices with diagonal scalar entries,
$(\cdot)\mT$ denotes transpose,
$(\cdot)^{-1}$ denotes inverse,
and ``$\triangleq$'' denotes equality by definition.
In addition, we write $\lambda_{\rom{min}}(A)$ (resp., $\lambda_{\rom{max}}(A)$) for the minimum (resp., maximum) eigenvalue of the Hermitian matrix $A$,
$\rom{tr}(\cdot)$ for the trace operator,
$\rom{vec}(\cdot)$ for the column stacking operator,
$\norm{\cdot}_2$ for the Euclidian norm,
$\norm{\cdot}_\infty$ for the infinity norm,
and $\norm{\cdot}_\rom{F}$ for the Frobenius matrix norm.
Furthermore, for a signal $x(t)=[x_1(t), x_2(t), \ldots, x_n(t)]\mT \in \IR^n$ defined for all $t \ge 0$,
the truncated $\mathcal{L}_\infty$ norm and the $\mathcal{L}_\infty$ norm [\citen{khalil}, Section 5] are defined as
$\norm{x_\tau(t)}_{\mathcal{L}_\infty} \triangleq \max_{1\le i \le n}(\sup_{0\le t \le \tau}|x_i(t)|)$
and $\norm{x(t)}_{\mathcal{L}_\infty} \triangleq \max_{1\le i \le n}(\sup_{t \ge 0}|x_i(t)|)$, respectively.


\ \section{Preliminaries}

Consider the nonlinear uncertain dynamical system given by
\bequ
	\dot{x}_\rom{p}(t)&=&A_\rom{p}x_\rom{p}(t)+B_\rom{p}\Lambda u(t) + B_\rom{p} \delta_\rom{p}\bigl(x_\rom{p}(t)\bigl), \quad x_\rom{p}(0)=x_\rom{p_0}, \label{unc:sys}
\eequ
where $x_\rom{p}(t)\in\IR^{n_\rom{p}}$ is the accessible state vector, $u(t)\in\IR^m$ is the control input, $\delta_\rom{p}:\IR^{n_\rom{p}}\rightarrow\IR^m$ is an \textit{uncertainty}, $A_\rom{p}\in\IR^{n_\rom{p} \times n_\rom{p}}$ is a known system matrix, $B_\rom{p}\in\IR^{n_\rom{p} \times m}$ is a known control input matrix, and $\Lambda \in \IR^{m \times m}_+\cap \mathbb{D}^{m \times m}$ is an \textit{unknown} control effectiveness matrix. 
Furthermore, we assume that the pair $(A_\rom{p},  B_\rom{p})$ is controllable and the uncertainty is parameterized as
\bequ
	\delta_\rom{p}\bigl(x_\rom{p}\bigl)&=&W_\rom{p}\mT \sigma_\rom{p}\bigl(x_\rom{p}\bigl), \quad x_\rom{p}\in \IR^{n_\rom{p}}, \label{unc:prm}
\eequ
where $W_\rom{p}\in \IR^{s\times m}$ is an \textit{unknown} weight matrix and $\sigma_\rom{p}:\IR^{n_\rom{p}} \rightarrow \IR^s$ is a known basis function\footnote{For the case where the basis function $\sigma_\rom{p}\bigl(x_\rom{p}\bigl)$ is \textit{unknown}, parameterization in (\ref{unc:prm}) can be relaxed, for example, by considering $\delta_\rom{p}\bigl(x_\rom{p}\bigl)=W_\rom{p}\mT \sigma_\rom{p}^\rom{nn}\bigl(V_\rom{p}\mT x_\rom{p}\bigl)+\varepsilon_\rom{p}^\rom{nn}\bigl(x_\rom{p}\bigl)$, $x_\rom{p}\in \mathcal{D}_{x_\rom{p}}$, where $W_\rom{p}\in \IR^{s\times m}$ and $V_\rom{p}\in \IR^{n_\rom{p}\times s}$ are \textit{unknown} weight matrices, $\sigma_\rom{p}^\rom{nn}:\mathcal{D}_{x_\rom{p}} \rightarrow \IR^s$ is a known basis composed of neural networks function approximators, $\varepsilon_\rom{p}^\rom{nn}:\mathcal{D}_{x_\rom{p}}\rightarrow\IR^m$ is an \textit{unknown} residual error, and $\mathcal{D}_{x_\rom{p}}$ is a compact subset of $\IR^{n_\rom{p}}$ \cite{lewis}.} 
of the form $\sigma_\rom{p}\bigl(x_\rom{p}\bigl)=\bigl[\sigma_\rom{p_1}\bigl(x_\rom{p}\bigl), \sigma_\rom{p_2}\bigl(x_\rom{p}\bigl),\ldots,  \sigma_\rom{p_s}\bigl(x_\rom{p}\bigl)\bigl]\mT$.

To address command following, let $c(t)\in\IR^{n_\rom{c}}$ be a given bounded piecewise continuous command and $x_\rom{c}(t)\in\IR^{n_\rom{c}}$ be the integrator state satisfying
\vspace{0.0cm}
\bequ
	\dot{x}_\rom{c}(t)&=&E_\rom{p}x_\rom{p}(t)-c(t), \quad x_\rom{c}(0)=x_\rom{c_0}, \label{int:state}
\eequ
where $E_\rom{p}\in\IR^{n_\rom{c} \times n_\rom{p}}$ allows to choose a subset of $x_\rom{p}(t)$ to be followed by $c(t)$\footnote{For stabilization, the integrator state given by (\ref{int:state}) may not be required (i.e., $n_\rom{c}=0$) since $c(t)\equiv0$ for this case.}. Now, (\ref{unc:sys}) can be augmented with (\ref{int:state}) as
\bequ
	\dot{x}(t)&\hspace{-0.1cm}=\hspace{-0.1cm}&Ax(t)+B\Lambda u(t)+BW_\rom{p}\mT\sigma_\rom{p}\bigl(x_\rom{p}(t)\bigl)+B_\rom{r}c(t), \quad x(0)=x_0, \label{unc:sys2}
\eequ
where $x(t)\triangleq\bigl[x_\rom{p}\mT(t),  x_\rom{c}\mT(t)\bigl]\mT\in\IR^n$, $n=n_\rom{p}+n_\rom{c}$, is the (augmented) state vector, $x_0\triangleq\bigl[x_\rom{p_0}\mT, x_\rom{c_0}\mT\bigl]\mT\in\IR^n$, 
\bequ
	A &\triangleq& \begin{bmatrix} A_\rom{p} & 0_{n_\rom{p}\times n_\rom{c}}  \\E_\rom{p} & 0_{n_\rom{c}\times n_\rom{c}}  \end{bmatrix} \in \IR^{n \times n}, \\
	B&\triangleq&\bigl[B_\rom{p}\mT, 0_{n_\rom{c}\times m}\mT\bigl]\mT  \in\IR^{n \times m}, \\
	B_\rom{r}&\triangleq&\bigl[0_{n_\rom{p}\times n_\rom{c}}\mT,-I_{n_\rom{c}\times n_\rom{c}} \bigl]\mT\in\IR^{n \times n_\rom{c}}.
\eequ

Next, consider the feedback control law given by
\bequ
	u(t) &=& u_\rom{n}(t)+u_\rom{a}(t), \label{total:ctrl}
\eequ
where $u_\rom{n} (t) \in \IR^m$ and $u_\rom{a}(t)\in\IR^m$ are the nominal and adaptive control laws, respectively. 
Furthermore, let the nominal control law be 
\bequ
	u_\rom{n}(t) &=& -Kx(t), \quad K\in\IR^{m \times n}, \label{nom:ctrl}
\eequ
such that $A_\rom{r} \triangleq A-BK$ is Hurwitz. 
Using (\ref{total:ctrl}) and (\ref{nom:ctrl}) in (\ref{unc:sys2}) yields
\bequ
	\dot{x}(t)&\hspace{-0.125cm}=\hspace{-0.125cm}&A_\rom{r}x(t)+B_\rom{r}c(t)+B\Lambda\bigl[u_\rom{a}(t)+W\mT\sigma\bigl(x(t)\bigl)\bigl], \ \ \ \ \label{unc:sys3}
\eequ
where $W\mT\triangleq\bigl[ \Lambda^{-1} W_\rom{p}\mT, (\Lambda^{-1}-I_{m\times m})K\bigl]\in\IR^{(s+n)\times m}$ is an \textit{unknown} (aggregated) weight matrix and $\sigma\mT\bigl(x(t)\bigl)\triangleq\bigl[\sigma_\rom{p}\mT\bigl(x_\rom{p}(t)\bigl), x\mT(t)\bigl]\in\IR^{s+n}$ is a known (aggregated) basis function. 
Considering (\ref{unc:sys3}), let the adaptive control law be
\bequ
	u_\rom{a}(t)=-\hat{W}\mT(t)\sigma\bigl(x(t)\bigl), \label{adapt:ctrl}
\eequ
where $\hat{W}(t)\in\IR^{(s+n)\times m}$ be the estimate of $W$ satisfying the update law
\bequ
	\dot{\hat{W}}(t)&=&\gamma \sigma\bigl(x(t)\bigl) e\mT(t) PB, \quad \hat{W}(0)=\hat{W}_0, \label{update:law} \ \ \
\eequ
where $\gamma \in \IR_+$ is the learning rate, $e(t)\triangleq x(t)-x_\rom{r}(t)$ is the system error with $x_\rom{r}(t)\in\IR^n$ being the reference state vector satisfying the reference system
\bequ
	\dot{x}_\rom{r}(t)&=&A_\rom{r}x_\rom{r}(t)+B_\rom{r}c(t), \quad x_\rom{r}(0)=x_\rom{r_0}, \label{ref:ideal}
\eequ
and $P \in \IR^{n \times n}_+ \cap \IS^{n \times n}$ is a solution of the Lyapunov equation\footnote{Since $A_\rom{r}$ is Hurwitz, it follows from [\citen{haddad}, Section 3.7] that there exists a unique $P$ satisfying (\ref{lyap:eqn}) for a given $R$.}
\vspace{-0.1cm}
\bequ
	0 &=& A_\rom{r}\mT P + P A_\rom{r} + R, \label{lyap:eqn}
\eequ
with $R \in \IR^{n \times n}_+ \cap \IS^{n \times n}$. 

Now, the system error dynamics is given by using (\ref{unc:sys3}), (\ref{adapt:ctrl}), and (\ref{ref:ideal}) as
\bequ
	\dot{e}(t)&\hspace{-0.05cm}=\hspace{-0.05cm}&A_\rom{r}e(t)-B\Lambda\tilde{W}\mT(t)\sigma\bigl(x(t)\bigl), \quad e(0)=e_0, \ \ \ \ \label{err:sys}
\eequ
where $\tilde{W}(t)\triangleq \hat{W}(t)-W \in \IR^{(s+n)\times m}$ is the weight error and $e_0 \triangleq x_0-x_\rom{r_0}$. 
The update law given by (\ref{update:law}) can be derived by using Lyapunov analysis by considering the Lyapunov function candidate
\bequ
	\mathcal{V}\bigl(e,\tilde{W}\bigl)&=& e\mT P e + \gamma^{-1}\rom{tr} \ \bigl(\tilde{W} \Lambda^\frac{1}{2}\bigl)\mT\bigl(\tilde{W} \Lambda^\frac{1}{2}\bigl). \ \ \label{lyap:cand}
\eequ
Note that $\mathcal{V}\bigl(0,0\bigl)=0$, $\mathcal{V}\bigl(e,\tilde{W}\bigl)>0$ for all $(e,\tilde{W})\neq(0,0)$, and 
$\mathcal{V}\bigl(e,\tilde{W}\bigl)$ is radially unbounded. Now, differentiating (\ref{lyap:cand}) and then using (\ref{update:law}) and (\ref{err:sys}) yields 
\bequ
\dot{\mathcal{V}}\bigl(e(t),\tilde{W}(t)\bigl)=-e\mT(t) R e(t)\le 0,
\eequ 
which guarantees that the system error $e(t)$ and the weight error $\tilde{W}(t)$ are Lyapunov stable, and hence, are bounded for all $t \in \overline{\IR}_+$. Since $\sigma\bigl(x(t)\bigl)$ is bounded for all $t \in \overline{\IR}_+$, it follows from (\ref{err:sys}) that $\dot{e}(t)$ is bounded, and hence, $\ddot{\mathcal{V}}\bigl(e(t),\tilde{W}(t)\bigl)$ is bounded for all $t \in \overline{\IR}_+$. Now, it follows from Barbalat's lemma [\citen{khalil}, Lemma 8.2] that 
\bequ
\lim_{t\rightarrow\infty}\dot{\mathcal{V}}\bigl(e(t),\tilde{W}(t)\bigl)=0,
\eequ 
which consequently shows that $\lim_{t\rightarrow\infty}e(t)=0$.

\textit{Remark 3.1}. Although $\lim_{t\rightarrow\infty}e(t)=0$, the state vector $x(t)$ can be far different 
from $x_\rom{r}(t)$ during transient time (learning phase), unless a high learning rate $\gamma$ is used in 
the update law (\ref{update:law}). 
To see this, we first let $x_\rom{r_0}=x_0$ in (\ref{ref:ideal}), and hence, $e(0)=0$ in (\ref{err:sys}). Note that this condition is realizable since it is assumed that the state vector $x_\rom{p}(t)$ in (\ref{unc:sys}) is accessible [\citen{initial:cond}, Section II.C]. Since $\dot{\mathcal{V}}\bigl(e(t),\tilde{W}(t)\bigl)\le 0$, and hence, 
\bequ
\mathcal{V}\bigl(e(t),\tilde{W}(t)\bigl)&\le& \mathcal{V}\bigl(e_0,\tilde{W}_0\bigl)=\gamma^{-1}\norm{\tilde{W}_0 \Lambda^\frac{1}{2}}^2_\rom{F}, \label{temp:01}
\eequ
then using $\mathcal{V}\bigl(e(t),\tilde{W}(t)\bigl)\ge\lambda_\rom{min}(P)\norm{e(t)}^2_2$ in (\ref{temp:01}) yields 
$\norm{e(t)}_2\le \norm{\tilde{W}_0 \Lambda^\frac{1}{2}}_\rom{F}/\sqrt{\gamma \lambda_\rom{min}(P)}$. 
Since $\norm{\cdot}_\infty \le \norm{\cdot}_2$ holds, and this bound is uniform, 
then $\norm{e_\tau(t)}_{\mathcal{L}_\infty}\le \norm{\tilde{W}_0 \Lambda^\frac{1}{2}}_\rom{F}/\sqrt{\gamma \lambda_\rom{min}(P)}$. 
Finally, since this holds uniformly in $\tau$, we have 
\bequ
\norm{e(t)}_{\mathcal{L}_\infty}\le \norm{\tilde{W}_0 \Lambda^\frac{1}{2}}_\rom{F}/\sqrt{\gamma \lambda_\rom{min}(P)},
\eequ 
which shows that the distance between $x(t)$ and $x_\rom{r}(t)$ can be made arbitrarily small in transient time by resorting to a high learning rate $\gamma$.
As discussed, however, update laws with high learning rates in the face of large system uncertainties and abrupt changes may yield to signals with high-frequency oscillations, which can violate actuator rate saturation constraints and/or excite unmodeled system dynamics [\citen{oscil:2}, \citen{oscil:3}] resulting in system instability for practical applications. 


\ \section{Frequency-Limited System Error Dynamics}

One of the fundamental components of a model reference adaptive control scheme is the system error $e(t)$. 
In particular, if the system error $e(t)$ contains any high-frequency oscillations, then the adaptive control law (\ref{adapt:ctrl}) can also have such oscillations, since the update law (\ref{update:law}) is driven by this system error $e(t)$. 
Motivating from this standpoint, our aim is to limit the frequency-content of the system error dynamics (\ref{err:sys}) during transient-time (learning phase), and hence, to filter out any possible high-frequency oscillations contained in the error signal $e(t)$. For this purpose, let $e_\rom{L}(t) \in \IR^n$ be a low-pass filtered system error of $e(t)$ given by
\bequ
	\dot{e}_\rom{L}(t)&\hspace{-0.1cm}=\hspace{-0.1cm}&A_\rom{r}e_\rom{L}(t)+\eta\bigl(e(t)-e_\rom{L}(t)\bigl), \quad e_\rom{L}(0)=0, \label{filter} \ \ \ \ \ \
\eequ
where $\eta\in\IR_+$ is a filter gain. 
Note that since $e_\rom{L}(t)$ is a low-pass filtered system error of $e(t)$, the filter gain $\eta$ is chosen such that $\eta\le\eta^*$, where $\eta^*\in\IR_+$ is a design parameter. 

Next, we add a mismatch term to the system error dynamics (\ref{err:sys}) in order to enforce a distance condition between the trajectories of the system error $e(t)$ and the trajectories of its low-pass filtered version $e_\rom{L}(t)$. 
This leads to a minimization problem involving an error criterion capturing the distance between $e(t)$ and $e_\rom{L}(t)$. In particular, consider the cost function given by
\vspace{-0.2cm}
\bequ
	\mathcal{J}\bigl(e,e_\rom{L}\bigl) &=& \frac{1}{2}\norm{e-e_\rom{L}}_2^2, \label{cost:func}
\eequ
and note that the negative gradient of (\ref{cost:func}) with respect to $e$ is given by
\vspace{-0.1cm}
\bequ
	\frac{\partial \bigl[ -\mathcal{J}\bigl(e(t),e_\rom{L}(t)\bigl) \bigl]}{\partial e(t)}=-\bigl(e(t)-e_\rom{L}(t)\bigl), \label{deriv:cost}
\eequ
which gives the structure of the proposed mismatch term. 
Using the idea presented in [\citen{whitaker:2},\citen{yucelen:1}--\citen{yucelen:2}], we now need to add (\ref{deriv:cost}) to the system error dynamics given by (\ref{err:sys}). 
For this purpose, we modify the reference system (\ref{ref:ideal}) as
\vspace{-0.1cm} 
\bequ
	\dot{x}_\rom{r}(t)&=&A_\rom{r}x_\rom{r}(t)+B_\rom{r}c(t)+\kappa\bigl(e(t)-e_\rom{L}(t)\bigl), \quad  x_\rom{r}(0)=x_\rom{r_0}, \label{ref:nonideal}
\eequ
where $\kappa\in\IR_+$, and hence, the system error dynamics is given by using (\ref{unc:sys3}), (\ref{adapt:ctrl}), and (\ref{ref:nonideal}) as
\bequ
	\dot{e}(t)&\hspace{-0.2cm}=\hspace{-0.2cm}&A_\rom{r}e(t)\hspace{-0.07cm}-\hspace{-0.07cm}B\Lambda\tilde{W}\mT(t)\sigma\bigl(x(t)\bigl)-\kappa\bigl(e(t)-e_\rom{L}(t)\bigl), \quad e(0)=e_0. \ \ \ \ \label{err:sys:mod}
\eequ
Finally, note for the rest of this paper that the update law (\ref{update:law}) is driven by the system error $e(t)=x(t)-x_\rom{r}(t)$, where $x_\rom{r}(t)$ is obtained from (\ref{ref:nonideal}) (not (\ref{ref:ideal})).

\textit{Remark 4.1}. The reference system (\ref{ref:nonideal}) captures a desired closed-loop dynamical system behavior modified by a mismatch term $\kappa\bigl(e(t)-e_\rom{L}(t)\bigl)$ representing the high-frequency content between the uncertain dynamical system and this reference system. 
Although this implies a modification of the ideal (unmodified) reference system (\ref{ref:ideal}) during transient time, as we see in the following sections, this mismatch term allows to limit the frequency content of the system error dynamics (\ref{err:sys:mod}), which is used to drive the adaptive controller. 
In other words, the purpose of our methodology is to prevent the update law from attempting to learn through the high-frequency content of the system error. 

\textit{Remark 4.2}. As it is noted, the filter gain $\eta$ needs to be chosen such that $\eta\le\eta^*$, 
where $\eta^*$ needs to be small enough to cut off the high-frequency content of $e(t)$. 
To see the negative effect of high filter gain, let $\eta$ be sufficiently large. 
Then, $e(t)-e_\rom{L}(t) \approx 0$ as a consequence of (\ref{filter}), 
and hence, we approximately recover the ideal (unmodified) reference system given by (\ref{ref:ideal}). 
In this case, the proposed approach converges to a standard model reference adaptive control scheme, which has practical limitations as discussed earlier in the presence of high learning rate $\gamma$. 
Furthermore, as a special case of $\eta=0$, the proposed approach converges to the approach documented in \cite{notable:3}, since $e_\rom{L}(t)\equiv0$ for all $t \in \overline{\IR}_+$ as a consequence of (\ref{filter}). 
Once again, as discussed earlier, this selection for the filter gain $\eta$ may result in poor transient performance in the presence of exogenous low-frequency persistent disturbances. 
Therefore, from a practical point of view, this imposes another constraint in the selection of filter gain such that it also needs to satisfy $\eta_*\le\eta$, where $\eta_*\in\IR_+$ needs to be large enough in order to suppress the effects of exogenous low-frequency persistent disturbances. This phenomenon is illustrated in the next remark for a special case as well as later in the paper for a more general case.

\begin{figure}[t!]  \center   \epsfig{file=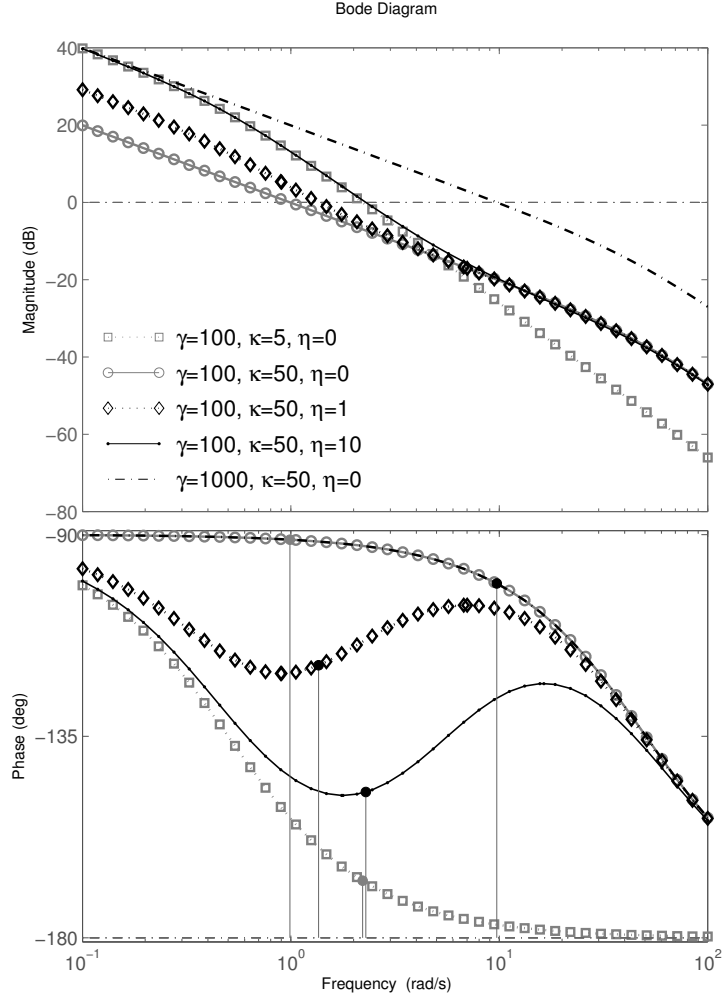, scale=0.6}
	\caption{Bode plots of the loop gain transfer function for different $\gamma$, $\kappa$, and $\eta$.}
	\label{figure:bode} \end{figure}

\textit{Remark 4.3}. To further elucidate the mechanism behind the proposed approach, 
let $n_\rom{p}=1$, $n_\rom{c}=0$, $m=1$, $A_\rom{p}=-\alpha$, $\alpha\in\IR_+$, $B_\rom{p}=\alpha$, $K=0$, 
and $\delta_\rom{p}\bigl(x_\rom{p}\bigl)=d$ with $d$ denoting an exogenous low-frequency disturbance. 
Furthermore, set $R=2$ in (\ref{lyap:eqn}) such that $P=\alpha^{-1}$ and let all initial conditions be zero. 
For this special case, the system loop transfer function $\mathcal{G}(s)$ (broken at the control input) can be equivalently written as a linear time-invariant dynamical system, and hence, we can resort to classical control theory tools, such as Bode plots, to analyze the closed-loop system with respect to different choices of $\gamma$, $\kappa$, and $\eta$. Specifically, the system loop transfer function is given by
\bequ
\mathcal{G}(s)&=&\underbrace{\frac{\gamma}{s}\biggl(\frac{s+\alpha+\eta}{s+\alpha+\kappa+\eta}\biggl)}_{\mathcal{C}(s)}\underbrace{\biggl(\frac{\alpha}{s+\alpha}\biggl)}_{\mathcal{P}(s)},
\eequ
where $\mathcal{C}(s)$ and $\mathcal{P}(s)$ denote the controller and the plant, respectively. 
Furthermore, for the standard model reference adaptive controller ($\eta$ is sufficiently large, or simply, $\kappa=0$), note that $\mathcal{C}(s)=\frac{\gamma}{s}$. 
For the controller $\mathcal{C}(s)$ of proposed approach, since $\frac{\gamma}{s}$ is multiplied by a lead compensator 
$\frac{s+\alpha+\eta}{s+\alpha+\kappa+\eta}$, it can improve stability margins of the closed-loop system positively. 
To further see the effects of $\gamma$, $\kappa$, and $\eta$, consider the Bode plot 
of $\mathcal{G}(s)$ in Figure \ref{figure:bode}.
Here, we tune these design parameters in order to obtain a large loop gain at low frequencies (from $0$ rad/s to $5/2\pi$ rad/s) for good rejection of low-frequency disturbances and a small loop gain at high frequencies to avoid injecting too much measurement noise into the plant [\citen{murray}, Section 11.4].
Furthermore, we also would like to have at least a time-delay margin of $0.25$ seconds. 
From Figure \ref{figure:bode}, one can see that the cases $(\gamma,\kappa,\eta)=(100, 5, 0)$, $(\gamma,\kappa,\eta)=(100, 50, 10)$, and $(\gamma,\kappa,\eta)=(1000, 50, 0)$ achieves approximately the same rejection of low-frequency disturbances. However, it should be noted that the case $(\gamma,\kappa,\eta)=(1000, 50, 0)$ amplifies measurement noise excessively in comparison to the cases $(\gamma,\kappa,\eta)=(100, 5, 0)$ and $(\gamma,\kappa,\eta)=(100, 50, 10)$. In addition, it should be also noted that the case $(\gamma,\kappa,\eta)=(100, 5, 0)$ has the poorest time-delay margin of $0.1$ seconds, whereas the cases $(\gamma,\kappa,\eta)=(100, 50, 10)$ and $(\gamma,\kappa,\eta)=(1000, 50, 0)$ have time-delay margins of $0.25$ and $0.14$ seconds, respectively. 
Therefore, one can conclude that the case $(\gamma,\kappa,\eta)=(100, 50, 10)$ achieves good rejection of low-frequency disturbances like the other two cases, has the maximum time-delay margin, and does not inject measurement noise as compared to the case $(\gamma,\kappa,\eta)=(1000, 50, 0)$. This shows the significance of having additional design parameters $\eta$ and $\kappa$ in the control design process. 
Moreover, the effect of increasing $\kappa$ alone can be depicted from the cases $(\gamma,\kappa,\eta)=(100, 5, 0)$ and $(\gamma,\kappa,\eta)=(100, 50, 0)$. 
Specifically, it deteriorates the rejection properties of low-frequency disturbances. 
That is the reason why we increased the adaptation gain $\gamma$ in the case $(\gamma,\kappa,\eta)=(1000, 50, 0)$ for achieving the same level of low-frequency disturbance rejection characteristics, however, as noted, this amplifies the measurement noise excessively and has less time-delay margin as compared to the case $(\gamma,\kappa,\eta)=(100, 50, 10)$. 
Finally, the effect of increasing $\eta$ to a moderate value can be seen from the cases $(\gamma,\kappa,\eta)=(100, 50, 0)$, $(\gamma,\kappa,\eta)=(100, 50, 1)$, and $(\gamma,\kappa,\eta)=(100, 50, 10)$. 
That is, we can recover the desired low-frequency disturbance rejection characteristics without increasing $\gamma$, and hence, without amplifying the measurement noise.


\ \section{Transient and Steady-State Performance Guarantees}

This section establishes transient and steady-state performance properties of the proposed adaptive control architecture. 
For this purpose, consider $e(t)=x(t)-x_\rom{r}(t)$ with $x_\rom{r}(t)$ satisfying (\ref{ref:nonideal}) and $\tilde{W}(t)=\hat{W}(t)-W$. 
Furthermore, let the ideal (unmodified) reference system\footnote{To prevent any abuse of notation, we redefine the ideal (unmodified) reference system in (\ref{ref:ideal}) as (\ref{ref:ideal:2}).} be
\bequ
	\dot{x}_\rom{r_i}(t)&\hspace{-0.05cm}=\hspace{-0.05cm}&A_\rom{r}x_\rom{r_i}(t)+B_\rom{r}c(t), \quad x_\rom{r_i}(0)=x_\rom{r_{0}}, \label{ref:ideal:2}
\eequ
where $x_\rom{r_i}(t)\in\IR^n$ being the ideal reference state vector. 
Finally, let $\tilde{x}(t)\triangleq x_\rom{r}(t)-x_\rom{r_i}(t)$ be the deviation error from the ideal reference system with $x_\rom{r}(t)$, once again, satisfying (\ref{ref:nonideal}). 
Then, the system error, weight update error, low-pass filtered system error, and the deviation error dynamics are, respectively, given by (\ref{err:sys:mod}), (\ref{filter}), 
\bequ
	\dot{\tilde{W}}(t)&=&\gamma \sigma\bigl(x(t)\bigl) e\mT(t) PB, \quad \tilde{W}(0)=\tilde{W}_0, \label{update:law:error} \\
	\dot{\tilde{x}}(t)&=&A_\rom{r} \tilde{x}(t)+\kappa\bigl(e(t)-e_\rom{L}(t)\bigl), \quad \tilde{x}(0)=0, \label{deviation:error} \ \ \ \
\eequ
where $\tilde{W}_0 \triangleq \hat{W}_0-W$.
The next theorem presents the first result of this paper.

\textit{Theorem 5.1}. 
Consider the nonlinear uncertain dynamical system given by (\ref{unc:sys}) subject to (\ref{unc:prm}), 
the (modified) reference system given by (\ref{ref:nonideal}), 
and the feedback control law given by (\ref{total:ctrl}) along with (\ref{nom:ctrl}), (\ref{adapt:ctrl}), and (\ref{update:law}). 
Then, the solution $\bigl(e(t),\tilde{W}(t),e_\rom{L}(t),\tilde{x}(t)\bigl)$ given by (\ref{err:sys:mod}), (\ref{filter}), (\ref{update:law:error}), and (\ref{deviation:error}) is Lyapunov stable for all 
$\bigl(e_0,\tilde{W}_0,0,0\bigl) \in \IR^n \times \IR^{(s+n)\times m}\times\IR^n\times\IR^n$ and $t \in \overline{\IR}_+$, 
and $\lim_{t\rightarrow\infty}\bigl(x(t)-x_\rom{r_i}(t)\bigl)=0$.
For $t \in \overline{\IR}_+$, in addition, 
\vspace{-0.2cm}
\bequ
	\norm{x(t)\hspace{-0.025cm}-\hspace{-0.025cm}x_\rom{r_i}(t)}_{\mathcal{L}_\infty}\hspace{-0.025cm}\le\hspace{-0.025cm}\sqrt{\frac{\epsilon_\mathcal{V}}{\lambda_\rom{min}(P)}}\Biggl(1\hspace{-0.1cm}+\hspace{-0.1cm}\sqrt{\frac{\kappa \lambda_\rom{max}(P)}{2\xi\lambda_\rom{min}(R)}}\Biggl), \label{bound:1}
\eequ
where $\xi\in(0,1)$ and $\epsilon_\mathcal{V}\triangleq\gamma^{-1}\norm{\tilde{W}_0\Lambda^\frac{1}{2}}_\rom{F}^2+\lambda_\rom{max}(P)\norm{e_0}_2^2$.

\textit{Proof}. To show Lyapunov stability of the solution $\bigl(e(t),\tilde{W}(t),e_\rom{L}(t),\tilde{x}(t)\bigl)$ given by (\ref{err:sys:mod}), (\ref{filter}), (\ref{update:law:error}), and (\ref{deviation:error}) for all 
$\bigl(e_0,\tilde{W}_0,0,0\bigl) \in \IR^n \times \IR^{(s+n)\times m}\times\IR^n\times\IR^n$ and $t \in \overline{\IR}_+$, 
consider the Lyapunov function candidate
\bequ
	\mathcal{V}^*\bigl(e,\tilde{W},e_\rom{L},\tilde{x}\bigl)&\hspace{-0.1cm}=\hspace{-0.1cm}&\mathcal{V}\bigl(e,\tilde{W}\bigl)+\eta^{-1}\kappa e_\rom{L}\mT P e_\rom{L}+2\xi \kappa^{-1} \lambda_\rom{max}^{-1}(P)\lambda_\rom{min}(R)\tilde{x}\mT P \tilde{x}, \label{proof:1} \ \ \ \ \
\eequ
where $\mathcal{V}\bigl(e,\tilde{W}\bigl)$ is given by (\ref{lyap:cand}), and note that 
$\mathcal{V}^*(0,0,0,0)=0$, $\mathcal{V}^*\bigl(e,\tilde{W},e_\rom{L},\tilde{x}\bigl)>0$ for all $\bigl(e,\tilde{W},e_\rom{L},\tilde{x}\bigl)$ $\neq(0,0,0,0)$, and $\mathcal{V}^*\bigl(e,\tilde{W},e_\rom{L},\tilde{x}\bigl)$ is radially unbounded. Differentiating (\ref{proof:1}) along the trajectories of (\ref{err:sys:mod}), (\ref{filter}), (\ref{update:law:error}), and (\ref{deviation:error}) yields
\bequ
	\dot{\mathcal{V}}^*(\cdot)&\hspace{-0.15cm}=\hspace{-0.15cm}&-e\mT(t)Re(t)-\eta^{-1}\kappa e_\rom{L}\mT(t)Re_\rom{L}(t)-2\xi \kappa^{-1}\lambda_\rom{max}^{-1}(P)\lambda_\rom{min}(R)\tilde{x}\mT(t) R \tilde{x}(t)\nonumber\\&&
	-2\kappa e\mT(t)Pe_\rom{H}(t) +2\kappa e_\rom{L}\mT(t)Pe_\rom{H}(t)+4\xi \lambda_\rom{max}^{-1}(P)\lambda_\rom{min}(R)\tilde{x}\mT(t) Pe_\rom{H}(t)\nonumber\\
	&\hspace{-0.15cm}=\hspace{-0.15cm}&-e\mT(t)Re(t)-\eta^{-1}\kappa e_\rom{L}\mT(t)Re_\rom{L}(t)-2\xi \kappa^{-1} \lambda_\rom{max}^{-1}(P)\lambda_\rom{min}(R)\tilde{x}\mT(t) R \tilde{x}(t)\nonumber\\&&
	-2\kappa e_\rom{H}\mT(t)Pe_\rom{H}(t) +2\xi\lambda_\rom{max}^{-1}(P)\lambda_\rom{min}(R)\Bigl[2\tilde{x}\mT(t)P^\frac{1}{2}P^\frac{1}{2}e_\rom{H}(t)  \Bigl]. \label{proof:3}
\eequ
where $e_\rom{H}(t)\triangleq e(t)-e_\rom{L}(t)$\footnote{Since $e_\rom{L}(t)$ is a low-pass filtered system error of $e(t)$, $e_\rom{H}(t)$ represents the high-frequency content of the system error.}.
Now, consider 
\bequ 
 2\tilde{x}\mT P^\frac{1}{2}P^\frac{1}{2}e_\rom{H}&\hspace{-0.1cm}\le\hspace{-0.1cm}&\mu \tilde{x}\mT P\tilde{x}+\frac{1}{\mu}e_\rom{H}\mT P e_\rom{H}, \quad \mu \in \IR_+, \ \ \ \ \label{proof:4}
 \eequ 
which follows from Young's inequality [\citen{berns}, Fact 1.4.7]. 
Using (\ref{proof:4}) in the last term of (\ref{proof:3}) and then etting $\mu=\xi\kappa^{-1}\lambda_\rom{max}^{-1}(P) \lambda_\rom{min}(R)$ in (\ref{proof:4}) yields
\bequ
	\dot{\mathcal{V}}^*(\cdot)&\hspace{-0.15cm}\le\hspace{-0.15cm}& -\lambda_\rom{min}(R)\norm{e(t)}_2^2
	- \eta^{-1}\kappa \lambda_\rom{min}(R)\norm{e_\rom{L}(t)}^2_2 \nonumber\\
	&& -2\xi\kappa^{-1} \lambda_\rom{min}^2(R)\lambda_\rom{max}^{-1}(P)\bigl[1-\xi\bigl]\norm{\tilde{x}(t)}^2_2, \label{proof:10}
\eequ
and hence, since $1-\xi>0$ in (\ref{proof:10}) by the definition of $\xi$, $\dot{\mathcal{V}}^*(\cdot)\le0$.
Therefore, the solution $\bigl(e(t),\tilde{W}(t),e_\rom{L}(t),\tilde{x}(t)\bigl)$ given by (\ref{err:sys:mod}), (\ref{filter}), (\ref{update:law:error}), and (\ref{deviation:error}) is Lyapunov stable for all 
$\bigl(e_0,\tilde{W}_0,0,0\bigl) \in \IR^n \times \IR^{(s+n)\times m}\times\IR^n\times\IR^n$ and $t \in \overline{\IR}_+$.

To show $\lim_{t\rightarrow\infty}\bigl(x(t)-x_\rom{r_i}(t)\bigl)=0$, note that $\sigma\bigl(x(t)\bigl)$ is bounded for all $t \in \overline{\IR}_+$, and hence, $\dot{e}(t)$ is bounded. Furthermore, since $\dot{e}_\rom{L}(t)$ and $\dot{\tilde{x}}(t)$ are also bounded, then $\ddot{\mathcal{V}}^*\bigl(e(t),\tilde{W}(t),e_\rom{L}(t),$ $\tilde{x}(t)\bigl)$ is bounded for all $t \in \overline{\IR}_+$. Now, it follows from Barbalat's lemma [\citen{khalil}, Lemma 8.2] that 
\bequ
\rom{lim}_{t \rightarrow \infty}\dot{\mathcal{V}}^*\bigl(e(t),\tilde{W}(t),e_\rom{L}(t),\tilde{x}(t)\bigl)=0,
\eequ 
which consequently shows that 
\bequ
\lim_{t\rightarrow\infty}e(t)&=&0, \\
\lim_{t\rightarrow\infty}e_\rom{L}(t)&=&0, \\
\lim_{t\rightarrow\infty}\tilde{x}(t)&=&0.
\eequ
Therefore, it follows from 
\bequ
x(t)-x_\rom{r_i}(t)&=&x(t)-x_\rom{r}(t)+x_\rom{r}(t)-x_\rom{r_i}(t)=e(t)+\tilde{x}(t),\label{in:ineq}
\eequ
that 
\bequ
\lim_{t\rightarrow\infty}\bigl(x(t)-x_\rom{r_i}(t)\bigl)&=&0.
\eequ

Finally, since $\dot{\mathcal{V}}^*\bigl(e(t),\tilde{W}(t),e_\rom{L}(t),\tilde{x}(t)\bigl)\le0$ for all $t \in \overline{\IR}_+$, 
this implies that 
\bequ
{\mathcal{V}}^*\bigl(e(t),\tilde{W}(t),e_\rom{L}(t),\tilde{x}(t)\bigl)&\le&{\mathcal{V}}^*\bigl(e_0,\tilde{W}_0,0,0\bigl). \ \ \ \label{proof:20}
\eequ
Using ${\mathcal{V}}^*\bigl(e_0,\tilde{W}_0,0,0\bigl)\le\epsilon_\mathcal{V}$ and
\bequ
{\mathcal{V}}^*\bigl(e(t),\tilde{W}(t),e_\rom{L}(t),\tilde{x}(t)\bigl)&\ge&\lambda_\rom{min}(P)\norm{e(t)}^2_2, \ \ \ 
\eequ
in (\ref{proof:20}) yields $\norm{e(t)}_2\le \sqrt{\lambda_\rom{min}^{-1}(P)\epsilon_\mathcal{V}}$. 
Since $\norm{\cdot}_\infty \le \norm{\cdot}_2$, and this bound is uniform, then 
$\norm{e_\tau(t)}_{\mathcal{L}_\infty}\le \sqrt{\lambda_\rom{min}^{-1}(P)\epsilon_\mathcal{V}}$, 
and hence, 
\bequ
\norm{e(t)}_{\mathcal{L}_\infty}\le \sqrt{\lambda_\rom{min}^{-1}(P)\epsilon_\mathcal{V}}, \label{temp:05}
\eequ 
is a direct consequence due to the fact that the former expression holds uniformly in $\tau$. 
Similarly, using ${\mathcal{V}}^*\bigl(e_0,\tilde{W}_0,0,0\bigl)\le\epsilon_\mathcal{V}$ and
\bequ
{\mathcal{V}}^*\bigl(e(t),\tilde{W}(t),e_\rom{L}(t),\tilde{x}(t)\bigl)&\ge&2\xi \kappa^{-1} \lambda_\rom{max}^{-1}(P) \lambda_\rom{min}(P) \lambda_\rom{min}(R)\norm{\tilde{x}(t)}^2_2,
\eequ
yields 
\bequ
\norm{\tilde{x}(t)}_{\mathcal{L}_\infty}\le\sqrt{   \frac{1}{2}\xi^{-1}\lambda_\rom{min}^{-1}(R) \lambda_\rom{min}^{-1}(P)  \kappa \lambda_\rom{max}(P)  \epsilon_\mathcal{V} }. \label{temp:06}
\eequ 
Now, it follows from (\ref{in:ineq}) that 
\bequ
\norm{x(t)-x_{\rom{r_i}}(t)}_{\mathcal{L}_\infty}\le\norm{e(t)}_{\mathcal{L}_\infty}+\norm{\tilde{x}(t)}_{\mathcal{L}_\infty}, \label{temp:07}
\eequ and hence, (\ref{bound:1}) is a direct consequence of using (\ref{temp:05}) and (\ref{temp:06}) in (\ref{temp:07}). This completes the proof. \hfill $\square$

\textit{Remark 5.1}. Even though the proposed architecture is predicated on a modified reference system given by (\ref{ref:nonideal}), Theorem 5.1 shows that $\lim_{t\rightarrow\infty}\bigl(x(t)-x_\rom{r_i}(t)\bigl)=0$, that is the (augmented) state vector $x(t)$ of (\ref{unc:sys2}) asymptotically converges to the ideal reference state vector $x_\rom{r_i}(t)$ of (\ref{ref:ideal:2}). Furthermore, during transient time (learning phase), the worst-case transient performance bound between $x(t)$ and $x_\rom{r_i}(t)$ is given by (\ref{bound:1}). 
To further elucidate this performance bound, we let $x_\rom{r_0}=x_0$ in (\ref{ref:nonideal}), and hence, $e(0)=0$ in (\ref{err:sys:mod}). As noted in Remark 3.1, since it is assumed that the state vector $x_\rom{p}(t)$ in (\ref{unc:sys}) is accessible, this condition is realizable [\citen{initial:cond}, Section II.C]. 
Now, denoting $\epsilon_{\mathcal{V}_1}\triangleq \norm{\tilde{W}_0\Lambda^\frac{1}{2}}_\rom{F}/\sqrt{\lambda_\rom{min}(P)}$ and $\epsilon_{\mathcal{V}_2}\triangleq \sqrt{  \frac{1}{2}\xi^{-1}\lambda_\rom{min}^{-1}(R)\lambda_\rom{min}(P)  }$, it follows from (\ref{bound:1}) that 
\bequ
	\norm{x(t)\hspace{-0.025cm}-\hspace{-0.025cm}x_\rom{r_i}(t)}_{\mathcal{L}_\infty}\hspace{-0.025cm}\le \gamma^{-\frac{1}{2}}\epsilon_{\mathcal{V}_1}\Bigl(1+\kappa^\frac{1}{2}\epsilon_{\mathcal{V}_2}\Bigl), \label{bound:elu}
\eequ
for all $t \in \overline{\IR}_+$. The performance bound in (\ref{bound:elu}) implies that the distance between $x(t)$ and $x_\rom{r_i}(t)$ can be made arbitrarily small in transient time by resorting to a high learning rate $\gamma$, similar to Remark 3.1 for the standard model reference adaptive control scheme. However, as we see in the next section, by increasing $\kappa$, we make the distance between $e(t)$ and $e_\rom{L}(t)$ sufficiently small in transient time, and hence, a high learning rate $\gamma$ subject to a high $\kappa$ does not yield to signals with high-frequency oscillations. Finally, it should be also noted from (\ref{bound:elu}) that keeping $\gamma$ constant but increasing $\kappa$ may result in a larger distance between $x(t)$ and $x_\rom{r_i}(t)$, and therefore, both should be increased simultaneously in order to keep this distance consistent during transient time. 


\ \section{Suppressing High-Frequency System Error Dynamics in Transient Time}

We now show that the high-frequency content of the system error $e_\rom{H}(t)=e(t)-e_\rom{L}(t)$ can be effectively suppressed as one increases $\kappa$ design parameter of the modified reference system (\ref{ref:nonideal}).
The next theorem presents the second result of this paper. 

\textit{Theorem 6.1}. Consider the system error dynamics given by (\ref{err:sys:mod}) and the low-pass filtered system error dynamics given by (\ref{filter}). 
Then,
\bequ
	e_\rom{H}(t,\kappa) &=& \rom{e}^{-\kappa t}e_{_\rom{H_0}}+\mathcal{O}(\kappa^{-1}), \label{sing:per1}
\eequ
holds for a sufficiently high $\kappa$, where $e_{_\rom{H_0}}\triangleq e(0)-e_\rom{L}(0)=e_0$\footnote{Recall that $e_\rom{L}(0)=0$ in (\ref{filter}).}. 

\textit{Proof}. Let $\varepsilon\triangleq\kappa^{-1}$. Then, (\ref{err:sys:mod}) and (\ref{filter}) can be equivalently written as
\bequ
	\varepsilon \dot{e}(t)&\hspace{-0.2cm}=\hspace{-0.2cm}&\varepsilon A_\rom{r}e(t)\hspace{-0.05cm}-\hspace{-0.05cm}\varepsilon B \Lambda \tilde{W}\mT(t)\sigma\bigl(x(t)\bigl)-\bigl(e(t)\hspace{-0.05cm}-\hspace{-0.05cm}e_\rom{L}(t)\bigl), \label{hf:1} \ \ \ \ \ \\
	\dot{e}_\rom{L}(t)&\hspace{-0.2cm}=\hspace{-0.2cm}& A_\rom{r} e_\rom{L}(t)+\eta\bigl(e(t)-e_\rom{L}(t)\bigl). \label{hf:2}
\eequ
Since setting $\varepsilon=0$ results in
\bequ
	0&=&e(t)-e_\rom{L}(t),  \label{hf:3} \\
	\dot{e}_\rom{L}(t)&=&A_\rom{r}e_\rom{L}(t), \label{hf:4}
\eequ
then the system given by (\ref{hf:1}) and (\ref{hf:2}) is said to be the singularly perturbed model form, 
where $e(t)=e_\rom{L}(t)$ captures the isolated root. 
To shift the quasi steady-state of $e(t)$ to the origin, consider 
\bequ
e_\rom{H}(t)=e(t)-e_\rom{L}(t),
\eequ 
as a change of variables, which yields 
\bequ
	\frac{\rom{d}e_\rom{H}(\tau)}{\rom{d} \tau}&=&-e_\rom{H}(\tau), \quad e_\rom{H}(0)=e_{_\rom{H_0}}, \label{adasd}
\eequ
where $\tau$ is related to the original $t$ through $\tau=t/\varepsilon$. 
Now, as a direct consequence of [\citen{khalil}, Theorem 11.2]
\bequ
	e(t,\varepsilon)&=&e^{A_\rom{r}t} e_\rom{L}(0)+ \rom{e}^{-t/\varepsilon}e_{_\rom{H_0}} +  \mathcal{O}(\varepsilon)=\rom{e}^{-t/\varepsilon}e_{_\rom{H_0}}+\mathcal{O}(\varepsilon), \\
	e_\rom{L}(t,\varepsilon)&=&e^{A_\rom{r}t} e_\rom{L}(0)+\mathcal{O}(\varepsilon)=\mathcal{O}(\varepsilon),
\eequ
which leads to the result given by (\ref{sing:per1}). This completes the proof. \hfill $\square$

\textit{Remark 6.1}. Note that (\ref{hf:4}) and (\ref{adasd}) are referred as the reduced-order system and the boundary-layer system, respectively, where the former describes the asymptotic behavior and the latter describes the transient behavior. 
In particular, Theorem 6.1 shows that the transient high-frequency content of the system error $e_\rom{H}(t)$ is globally exponentially stable for a sufficiently high $\kappa$, and hence, it vanishes in a fast manner. 
If, in addition, $x_\rom{r_0}=x_0$ in (\ref{ref:nonideal}), then it follows from (\ref{sing:per1}) that 
$e_\rom{H}(t,\kappa) = \mathcal{O}(\kappa^{-1})$. 


\ \section{Extensions to Time-Varying Uncertainties}
This section discusses robustness properties of the proposed adaptive control architecture with respect to time-varying uncertainties. 
For this purpose, we relax the uncertainty parametrization given by (\ref{unc:prm}) to 
\bequ
	\delta_\rom{p}\bigl(t, x_\rom{p}\bigl)&=&W_\rom{p}\mT(t) \sigma_\rom{p}\bigl(x_\rom{p}\bigl), \quad x_\rom{p}\in \IR^{n_\rom{p}}, \label{unc:prm:rlx}
\eequ
where $W_\rom{p}(t)\in \IR^{s\times m}$ is an \textit{unknown} time-varying weight matrix subject to $\norm{W_\rom{p}(t)}_\rom{F}\le w_\rom{p,max}$, $w_\rom{p,max} \in \IR_+$, and $\norm{\dot{W}_\rom{p}(t)}_\rom{F}\le \dot{w}_\rom{p.max}$, $\dot{w}_\rom{p,max}\in\IR_+$ \footnote{If we let the first entry of the basis function $\sigma_\rom{p}\bigl(x_\rom{p}\bigl)$ to be the bias term, then this parameterization also captures the effect of exogenous time-varying disturbances.}. 
Note that, as a consequence of (\ref{unc:prm:rlx}), the \textit{unknown} (aggregated) weight matrix has the form 
$W(t)\mT\triangleq\bigl[ \Lambda^{-1} W_\rom{p}\mT(t), (\Lambda^{-1}-I_{m\times m})K\bigl]\in\IR^{(s+n)\times m}$ for this case. We introduce the following definition for developing the main results of this section. 

\textit{Definition 7.1}.
Let $\phi:\IR^{n } \rightarrow \IR$ be a continuously differentiable convex function given by
\bequ
\phi(\theta) \triangleq \frac{(\varepsilon_\theta+1)\theta\mT \theta-\theta^2_\rom{max}}{\varepsilon_\theta \theta^2_\rom{max}},
\eequ
where $\theta_\rom{max} \in \IR$ is a projection norm bound imposed on $\theta \in \IR^n$ and
$\varepsilon_\theta >0$ is a projection tolerance bound.
Then, the \textit{projection operator} $\rom{Proj}:\IR^{n } \times \IR^{n } \rightarrow \IR^{n }$ is defined by
\vspace{-0.1cm}
\begin{eqnarray}
   \Proj(\theta,y)
   \teq
   \left\{ \begin{array}{cl}
      \hspace{-2.8cm} y,  \ \ \ \mbox{ if $\phi(\theta) < 0$},
      \\
      \hspace{ 0.0cm} y, \ \ \ \mbox{ if $\phi(\theta) \ge 0$ and $\phi'(\theta)y  \le 0$},
      \\
      \hspace{-2.4cm} y - \frac{\phi'^\mathrm{^\rom{T}}(\theta)\phi'(\theta) y} {\phi'(\theta)\phi'^\mathrm{^\rom{T}}(\theta)} \phi(\theta), \\ \hspace{0.85cm} \mbox{ if $\phi(\theta) \ge 0$ and $\phi'(\theta)y>0 $},
   \end{array} \right.
   \label{my_proj_operator}
\end{eqnarray}
where $y \in \IR^{n }$.

\textit{Remark 7.1}. It follows from Definition 7.1 that 
\bequ
(\theta-\theta^*)\mT(\rom{Proj}(\theta,y)-y)\le0, \quad  \theta^*\in\IR^n, \label{defnt:b}
\eequ
holds \cite{project}.
The definition of the projection operator can be generalized to matrices as
\bequ
\rom{Proj}_\rom{m}(\Theta,Y)=\bigl( \rom{Proj}(\rom{col}_1 (\Theta), \rom{col}_1 (Y)), \ \ldots, \  \rom{Proj}(\rom{col}_m (\Theta), \rom{col}_m (Y)) \bigl),
\eequ
where $\Theta \in \IR^{n \times m}$, $Y \in \IR^{n \times m}$, and $\rom{col}_i (\cdot)$ denotes $i$-th column operator.
In this case, for a given $\Theta^* \in \IR^{n \times m}$, it follows from (\ref{defnt:b}) that
\bequ
\ \ \rom{tr} \bigl[ (\Theta\hspace{-0.1cm}-\hspace{-0.1cm}\Theta^*)\mT(\rom{Proj}_\rom{m}(\Theta,Y)\hspace{-0.1cm}-\hspace{-0.1cm}Y) \bigl ]=\sum_{i=1}^m \bigl[ \rom{col}_i(  \Theta \hspace{-0.1cm}-\hspace{-0.1cm}  \Theta^* )\mT \bigl(\rom{Proj}( \rom{col}_i (\Theta), \rom{col}_i (Y) )\hspace{-0.1cm} -\hspace{-0.1cm} \rom{col}_i (Y)   \bigl) \bigl]\le 0, \nonumber\\ \eequ
holds.
Throughout this section, we assume that the projection norm bound imposed on each column of $\Theta \in \IR^{n \times m}$ is $\theta_\rom{max}$.

Next, let the adaptive control law be given by (\ref{adapt:ctrl}), where $\hat{W}(t)\in\IR^{(s+n)\times m}$ is the estimate of $W(t)$ satisfying the update law
\bequ
	\dot{\hat{W}}(t)&=&\gamma \rom{Proj}_\rom{m} \Bigl[\hat{W}(t), \sigma\bigl(x(t)\bigl) e\mT(t) PB\Bigl],   \quad \hat{W}(0)=\hat{W}_0, \label{update:law:2} \ \ \
\eequ
with $\gamma \in \IR_+$ being the learning rate, $e(t)\triangleq x(t)-x_\rom{r}(t)$ being the system error such that $x_\rom{r}(t)\in\IR^n$ satisfies (\ref{ref:nonideal}), and $P \in \IR^{n \times n}_+ \cap \IS^{n \times n}$ being a solution of (\ref{lyap:eqn}). 

We note here that the system error, weight update error, low-pass filtered system error, and the deviation error dynamics are, respectively, given by (\ref{err:sys:mod}), (\ref{filter}), 
\bequ
	\dot{\tilde{W}}(t)&=&\gamma \rom{Proj}_\rom{m} \Bigl[\hat{W}(t), \sigma\bigl(x(t)\bigl) e\mT(t) PB\Bigl]-\dot{W}(t), \quad \tilde{W}(0)=\tilde{W}_0, \label{update:law:err2} \ \ \
\eequ
and (\ref{deviation:error}), where $\tilde{W}_0 \triangleq \hat{W}_0-W$.
Furthermore, since $W_\rom{p}(t)$ and $\dot{W}_\rom{p}(t)$ are bounded, there exists norm bounds $w_\rom{max}$ and $\dot{w}_\rom{max}$ such that $\norm{W(t)}_\rom{F}\le w_\rom{max}$ and $\norm{\dot{W}(t)}_\rom{F}\le \dot{w}_\rom{max}$ for all $t \in \overline{\IR}_+$. 
Similarly, since the update law for $\hat{W}(t)$ is predicated on the projection operator and $W(t)$ is bounded, 
there also exists a norm bound $\tilde{w}_\rom{max}$ such that $\norm{\tilde{W}(t)}_\rom{F}\le\tilde{w}_\rom{max}$ for all $t \in \overline{\IR}_+$. 
Finally, as it is standard for the projection operator-based update laws, we assume that $\norm{\hat{W}(0)}_\rom{F}\le \hat{w}_\rom{max}$, where $\hat{w}_\rom{max}$ is the projection norm bound imposed on $\hat{W}(t)$.
The next theorem presents the third result of this paper.

\textit{Theorem 7.1}.
Consider the nonlinear uncertain dynamical system given by (\ref{unc:sys}) subject to (\ref{unc:prm:rlx}), 
the (modified) reference system given by (\ref{ref:nonideal}), 
and the feedback control law given by (\ref{total:ctrl}) along with (\ref{nom:ctrl}), (\ref{adapt:ctrl}), and (\ref{update:law:2}). 
Then, the solution $\bigl(e(t),\tilde{W}(t),e_\rom{L}(t),\tilde{x}(t)\bigl)$ given by (\ref{err:sys:mod}), (\ref{filter}), (\ref{update:law:err2}), and (\ref{deviation:error}) is uniformly bounded for all 
$\bigl(e_0,\tilde{W}_0,0,0\bigl) \in \IR^n \times \IR^{(s+n)\times m}\times\IR^n\times\IR^n$ and $t \in \overline{\IR}_+$, 
with ultimate bound
\vspace{-0.1cm}
\bequ
	\norm{x(t)\hspace{-0.025cm}-\hspace{-0.025cm}x_\rom{r_i}(t)}_{\mathcal{L}_\infty}\hspace{-0.025cm}\le\hspace{-0.025cm}\sqrt{\frac{\rho_\mathcal{V}}{\lambda_\rom{min}(P)}}\Biggl(1\hspace{-0.1cm}+\hspace{-0.1cm}\sqrt{\frac{\kappa \lambda_\rom{max}(P)}{2\xi\lambda_\rom{min}(R)}}\Biggl), \label{ultim:bound}
\eequ
where $\xi\in(0,1)$ and 
\bequ
\rho_\mathcal{V}&\triangleq&\gamma^{-1}\norm{\Lambda}_\rom{F}\tilde{w}^2_\rom{max}\Biggl( 1+4\gamma^{-1}\norm{\Lambda}_\rom{F}\dot{w}^2_\rom{max}\lambda_\rom{min}^{-2}(R) \biggl[ 1+ \eta\kappa^{-1}\lambda_\rom{max}(P)\nonumber\\&&+\frac{1}{2}\kappa \lambda_\rom{max}^2(P) \lambda_\rom{min}(R)\xi^{-1}(1-\xi)^{-2}    \biggl] \Biggl).
\eequ

\textit{Proof.} Consider the Lyapunov-like function candidate given by (\ref{proof:1}). 
Then, using the property of the projection operator given by (\ref{defnt:b}) and following similar steps in the proof of Theorem 5.1, we have
\bequ
	\dot{\mathcal{V}}^*\bigl(\cdot\bigl)&\hspace{-0.2cm}\le\hspace{-0.2cm}&-c_1\norm{e(t)}_2^2-c_2\norm{e_\rom{L}(t)}^2_2-c_3\norm{\tilde{x}(t)}_2^2+c_4, \ \ \ \ \ \label{pr:1}
\eequ
where $c_1\triangleq\lambda_\rom{min}(R)$, $c_2 \triangleq \eta^{-1}\kappa \lambda_\rom{min}(R)$, 
$c_3 \triangleq 2 \xi \kappa^{-1}$ $\cdot \lambda_\rom{min}^2(R)\lambda_\rom{max}^{-1}(P)(1-\xi)$, 
and $c_4 \triangleq 2 \gamma^{-1} \norm{\Lambda}_\rom{F} \tilde{w}_\rom{max} \dot{w}_\rom{max}$, 
and hence, either
$\norm{e(t)}_2 \ge \psi_1 \triangleq \bigl(c_4/c_1\bigl)^2$ or 
$\norm{e_\rom{L}(t)}_2 \ge \psi_2 \triangleq \bigl(c_4/c_2\bigl)^2$ or 
$\norm{\tilde{x}(t)}_2 \ge \psi_3 \triangleq \bigl(c_4/c_3\bigl)^2$ yields $\dot{\mathcal{V}}^*(\cdot)\le0$. 
Since the update law for $\hat{W}(t)$ is predicated on the projection operator, $\norm{\tilde{W}(t)}\le\tilde{w}_\rom{max}$. 
That is, $\dot{\mathcal{V}}^*(\cdot)\le0$ outside the compact set defined by
\bequ
	\mathcal{D}_{\mathcal{V}} &\triangleq& \{(e(t),\tilde{W}(t),e_\rom{L}(t),\tilde{x}(t))\in\IR^n \times \IR^{(s+n) \times m}\times\IR^n\times\IR^n: \norm{e(t)}_2\le\psi_1 \nonumber\\
	&& \mathrm{or} \ \norm{e_\rom{L}(t)}_2\le\psi_2   \mathrm{or} \ \norm{\tilde{x}(t)}_2\le\psi_3\} \ \union \ \{(e(t),\tilde{W}(t),e_\rom{L}(t),\tilde{x}(t))\in\IR^n \nonumber\\&& \times \IR^{(s+n) \times m}\times\IR^n\times\IR^n: \norm{\tilde{W}(t)} \le\tilde{w}_\rom{max}\}. \label{set:set}
\eequ
Therefore, $\mathcal{V}^*\bigl(e(t),\tilde{W}(t),e_\rom{L}(t),\tilde{x}(t)\bigl)$ cannot grow outside $\mathcal{D}_\mathcal{V}$, and hence, 
\bequ
\mathcal{V}^*\bigl(e(t),\tilde{W}(t),e_\rom{L}(t),\tilde{x}(t)\bigl)&\le& \rom{max}_{(e,\tilde{W},e_\rom{L},\tilde{x})\in\mathcal{D}_\mathcal{V}}(e,\tilde{W},e_\rom{L},\tilde{x}\bigl)=\rho_\mathcal{V}.
\eequ 
Now, by following the final steps in the proof of Theorem 5.1, (\ref{ultim:bound}) is immediate. 
This completes the proof. \hfill $\square$

\textit{Remark 7.2}. By following an identical analysis to [\citen{notable:4}, Theorem 9], the results of Theorem 7.1 can be further extended to show that the Lyapunov-like function candidate given by (\ref{proof:1}) exponentially converges to the set defined by the ultimate bound (\ref{ultim:bound}).


\section{Illustrative Numerical Example}

\begin{figure}[] \center \epsfig{file=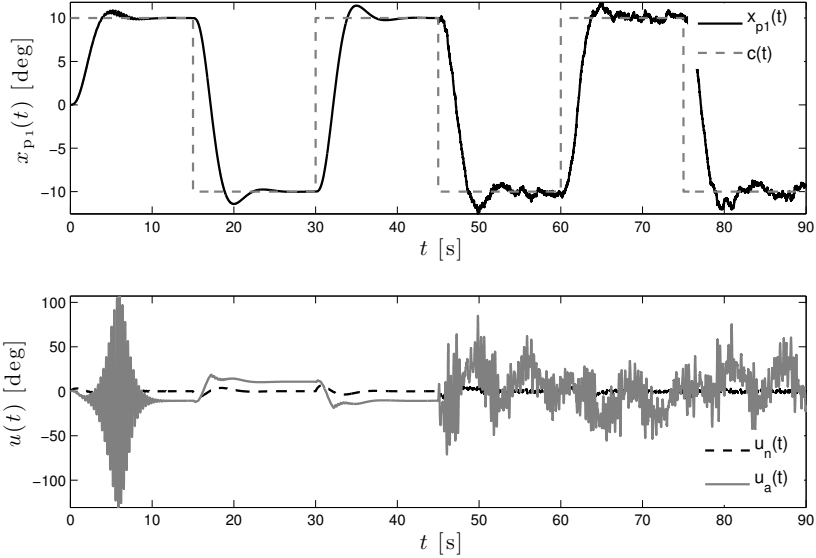, scale=0.56}
\caption{Command following performance for the standard model reference adaptive control approach ($\gamma=500$, $\kappa=0$, and $\eta=0$).}
\label{figure:sim:1} \end{figure}

\begin{figure}[] \center \epsfig{file=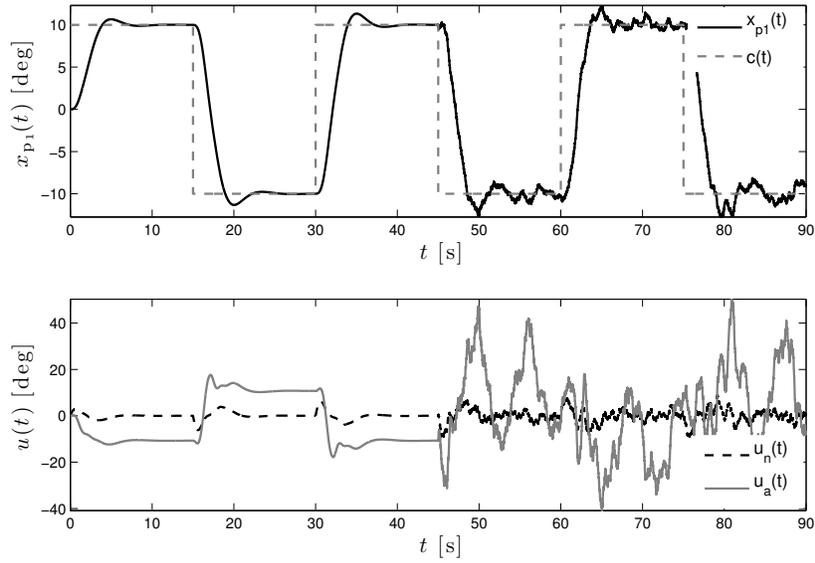, scale=0.56}
\caption{Command following performance for the proposed model reference adaptive control approach ($\gamma=500$, $\kappa=100$, and $\eta=5$).}
\label{figure:sim:2} \end{figure}

\begin{figure}[]   \center \epsfig{file=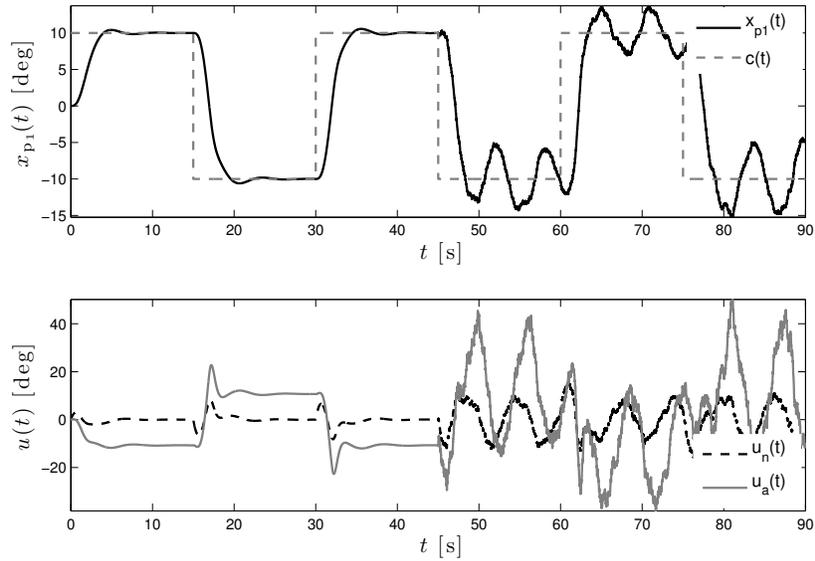, scale=0.56}
\caption{Command following performance for the model reference adaptive control approach of [\citen{notable:3},\citen{notable:4}] ($\gamma=500$, $\kappa=100$, and $\eta=0$).}
\label{figure:sim:3} \end{figure}

\begin{figure}[t!]  \center \epsfig{file=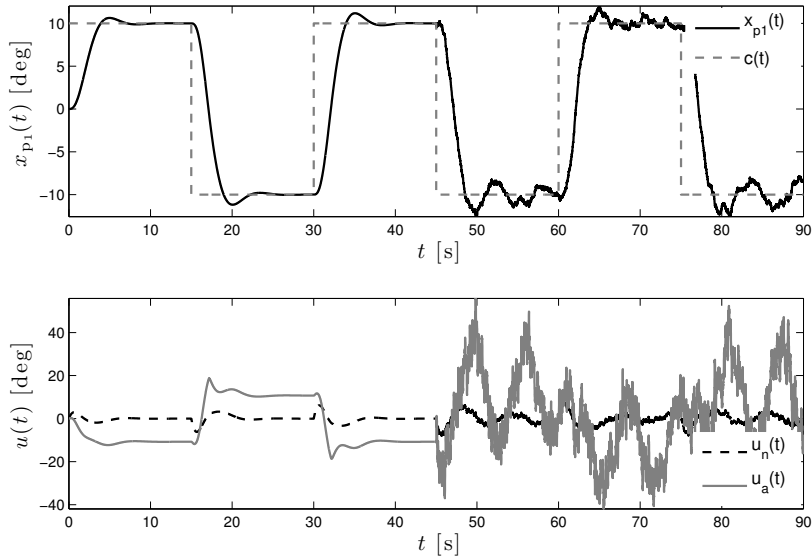, scale=0.56}
\caption{Command following performance for the model reference adaptive control approach of [\citen{notable:3},\citen{notable:4}] ($\gamma=2000$, $\kappa=100$, and $\eta=0$).}
\label{figure:sim:4} \end{figure}

Consider the nonlinear dynamical system representing a controlled wing rock aircraft dynamics model given by
\bequ
    \begin{bmatrix} \dot{x}_\rom{p_1}(t) \\ \dot{x}_\rom{p_2}(t) \end{bmatrix} &\hspace{-0.02cm}=\hspace{-0.02cm}& \begin{bmatrix} 0 & 1  \\0 & 0  \end{bmatrix} \begin{bmatrix} {x}_\rom{p_1}(t) \\ {x}_\rom{p_2}(t) \end{bmatrix} \hspace{-0.1cm}+\hspace{-0.1cm} \begin{bmatrix} 0 \\ 1 \end{bmatrix} \bigl[\Lambda u(t)+\delta_\rom{p}(t,x_\rom{p}(t))\bigl],\label{083}
\eequ
where $x_\rom{p_1}(0)=0$, $x_\rom{p_2}(0)=0$, $x_\rom{p_1}$ represents the roll angle in radians, and $x_\rom{p_2}$ represents the roll rate in radians per second.
In (\ref{083}), $\delta_\rom{p}(t,x_\rom{p})$ and $\Lambda$ represent uncertainties of the form
\bequ
\delta_\rom{p}(t,x_\rom{p})=\alpha_1\rom{sin}(t)+\alpha_2x_\rom{p_1}+\alpha_3x_\rom{p_2}+\alpha_4|{x_\rom{p_1}}|x_\rom{p_2}+\alpha_5|{x_\rom{p_2}}| x_\rom{p_2}+\alpha_6x_\rom{p_1}^3,
\eequ 
and 
\vspace{-0.25cm}
\bequ
\Lambda=0.75, 
\eequ 
where $\alpha_i$, $i=1,\ldots,6$, are unknown parameters.
For our numerical example, we set $\alpha_1=0.25$, $\alpha_2=0.5$, $\alpha_3=1.0$, $\alpha_4=-5.0$, $\alpha_5=5.0$, and $\alpha_6=10.0$. 
We chose $K=[2.0, \ 2.0, \ 1.0]$ for the nominal controller design. 
For the proposed adaptive control architecture (Theorem 7.1), 
$
\sigma(x)=\bigl[1, \ x_\rom{p_1}, \ x_\rom{p_2}, \ |{x_\rom{p_1}}|x_\rom{p_2}, \ |{x_\rom{p_2}}| x_\rom{p_2}, \ x_\rom{p_1}^3, \ x\mT \bigl]\mT,
$
is chosen as the basis function and we set $R=I_3$. 
Figures \ref{figure:sim:1}--\ref{figure:sim:4} present the results, where measurement noise is added to the state vector of (\ref{083}) and $\alpha_1$ is set from 0 to 0.25 both at $t=45$ seconds\footnote{That is, exogenous time-varying disturbance $\rom{sin}(t)$ is added at $t=45$ seconds.} for all cases. 
Here, our aim is to follow a given square-wave roll angle command $c(t)$. 

Figure \ref{figure:sim:1} shows the closed-loop system performance of the standard model reference adaptive control approach ($\gamma=500$, $\kappa=0$, and $\eta=0$). 
Even though we achieve a satisfactory command following performance with this approach, as discussed in Remark 3.1, its control performance is unacceptable due to high-frequency oscillations and measurement noise amplification.

Next, we show the closed-loop system performance of the proposed model reference adaptive control approach ($\gamma=500$, $\kappa=100$, and $\eta=5$) in Figure \ref{figure:sim:2}. 
In particular, we achieve a satisfactory command following performance similar to the case in Figure \ref{figure:sim:1}. 
However, the control response of our approach is clearly superior as compared to the control response of the standard model reference adaptive control approach in Figure \ref{figure:sim:1}. 
This is expected from the proposed theory, and hence, the control response of the proposed approach neither has high-frequency oscillations nor high measurement noise amplification.

In order to compare our approach with the approach of [\citen{notable:3},\citen{notable:4}], we set $\eta=0$ in Figure \ref{figure:sim:3}.  In this case, however, the transient performance is not sufficient due to the presence of exogenous low-frequency persistent disturbance. In order to improve the transients of this approach, we increased learning rate to $\gamma=2000$ in Figure \ref{figure:sim:4}. Even though we now achieve a satisfactory command following performance with this approach, its control response has excessive measurement noise as compared to the control response of the proposed approach in Figure \ref{figure:sim:2}.  
Note that this is also observed in Remark 4.3 for the case of a simple example.

\vspace{0.5cm}
\ \section{Conclusion}

We contributed to the previous studies in model reference adaptive control theory by introducing a new reference system in order to improve the transient performance. 
By utilizing singular perturbation theory, it is shown that the proposed reference system allows to limit the frequency content of the system error dynamics to yield fast adaptation without incurring high-frequency oscillations in the transient performance. 
We derived the guaranteed performance bounds, analyzed the effects of design parameters on the system performance, and discussed robustness properties to time-varying uncertainties and disturbances.


\ \section{Acknowledgment}

The authors wish to thank Eugene Lavretsky from the Boeing Company for his constructive suggestions. 

\vspace{-2cm}

\bibliographystyle{IEEEtran} 

\end{document}